# The Alpha-Beta-Skew-Logistic Distribution And Its Applications

Sricharan Shah and Partha Jyoti Hazarika

Department of Statistics, Dibrugarh University, Dibrugarh, Assam, 786004

**Abstract**

In this paper, an alpha-beta-skew-logistic distribution is proposed following the same methodology as those of alpha-beta-skew-normal of Shafiei et al. (2016) and investigated some of its related distributional properties. Finally, the validity of our proposed distribution has tested by considering three real life applications and comparing the values of Akaike Information Criterion (AIC) and Bayesian Information Criterion (BIC) with the values of some other related distributions. Likelihood ratio test is used for discriminating between logistic and the proposed distributions.

## 1. Introduction

Azzalini skew normal distribution was first pioneered by Azzalini (1985) with the following density function, denoted by $Z \sim SN(\lambda)$ and is given by

$$f_Z(z;\lambda) = 2\varphi(z)\Phi(\lambda z); \quad -\infty < z < \infty, -\infty < \lambda < \infty \qquad (1)$$

Where, $\varphi$ and $\Phi$ are the pdf and cdf of $N(0,1)$, respectively, and $\lambda$, the skewness parameter. In discussant with Arnold and Beaver (2002), Balakrishnan (2002) proposed the generalization of eqn. (1) and its density function is given by

$$f_Z(z;\lambda,n) = \varphi(z)[\Phi(\lambda z)]^n / C_n(\lambda) \qquad (2)$$

where, $n$ is positive integer and $C_n(\lambda) = E[\Phi^n(\lambda U)]$, $U \sim N(0,1)$. This distribution and its generalizations were studied by different researchers (for details see see Chakraborty and Hazarika, 2011).

For the construction of skew-symmetric distributions, Huang and Chen (2007) proposed a concept of skew function $G(.)$, a Lebesgue measurable function such that, $0 \leq G(z) \leq 1$ and $G(z) + G(-z) = 1$, $z \in R$, almost everywhere. Its density function is given by

$$f(z) = 2\varphi(z)G(z) \; ; z \in R \qquad (3)$$

By selecting an appropriate skew function $G(.)$ in the eqn. (3), one can construct different skew distributions. Therefore, considering the idea of eqn. (3), Elal-Olivero (2010) and Shafiei et al. (2016) developed alpha skew normal distribution and alpha beta skew normal distribution.

In this study the main aim is to develop a new distribution known as an alpha-beta-skew-logistic distribution following the idea of Balakrishnan (2002) and Shafiei et al. (2016) and investigate some of its basic properties. The three real life applications are used to this new proposed distribution which shows better fitting as compared to the other known related distributions considered here in this study. Therefore the remainder of the paper is organized as follows: In section 2 the density and distribution function are given. The moments, moment generating function and characteristics function are described in section 3. The skewness and kurtosis are given in section 4 with symmetric component random variable of the model. The estimation and applications of the proposed distribution with the other distributions are discussed in section 5 with likelihood ratio test. Conclusion is given in section 6.

## 2. Density and Distribution functions of AbSLG ($\alpha, \beta$):

**Definition 1:** A random variable $Z \sim AbSLG(\alpha, \beta)$ distribution if its density function defined as

$$f_Z(z; \alpha, \beta) = \frac{[(1 - \alpha z - \beta z^3)^2 + 1]}{C(\alpha, \beta)} \frac{e^{-z}}{(1 + e^{-z})^2} \; ; \; z, \alpha, \beta \in R \qquad (4)$$

where, $C(\alpha, \beta) = (210 + 35\pi^2 \alpha^2 + 98\pi^4 \alpha \beta + 155\pi^6 \beta^2)/105$. Then, it is said to be alpha-beta skew logistic distribution with skewness parameter $\alpha$ and $\beta$. In the rest of this article we shall refer the distribution in (4) as $AbSLG(\alpha, \beta)$.

**Some properties of AbSLG ($\alpha, \beta$):**

(i) If $\beta = 0$, then we get the alpha skew logistic distribution as

$$f_Z(z; \alpha) = \frac{105 \, e^{-z}[(1 - \alpha z)^2 + 1]}{(1 + e^{-z})^2 (210 + 35\pi^2 \alpha^2)}$$

(ii) If $\alpha = 0$, then we get $f_Z(z; \beta) = \frac{105 \, e^{-z}[(1 - \beta z^3)^2 + 1]}{(1 + e^{-z})^2 (210 + 155 \pi^6 \beta^2)}$

This equation is known as beta skew logistic ($bSLG(\beta)$) distribution.

(iii) If $\alpha = \beta = 0$, then we get the standard logistic distribution and is given by

$$f_Z(z) = \frac{e^{-z}}{(1+e^{-z})^2}$$

(iv) If $\alpha \to \pm\infty$, then we get the Bimodal Logistic ($BLG(2)$) distribution as

$$f_Z(z) = \frac{3z^2}{\pi^2} \frac{e^{-z}}{(1+e^{-z})^2}$$

(v) If $\beta \to \pm\infty$, then we get $f_Z(z) = \frac{21 z^6}{31 \pi^6} \frac{e^{-z}}{(1+e^{-z})^2}$

This equation is known as Bimodal Logistic ($BLG(6)$) distribution.

(vi) If $Z \sim AbSLG(\alpha,\beta)$, then $-Z \sim AbSLG(-\alpha,-\beta)$.

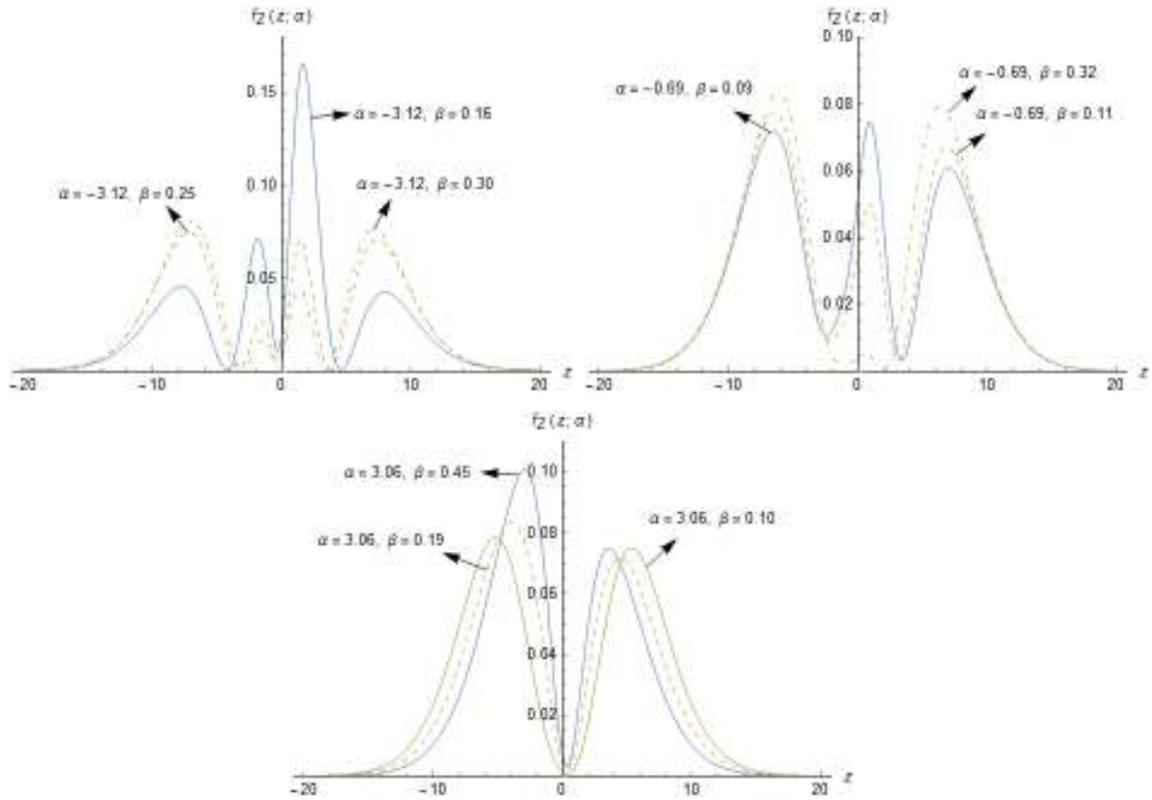

**Figure 1:** The plots of pdf of $AbSLG_2(\alpha,\beta)$ for different values of $\alpha$ and $\beta$

**Proposition 1:** The cumulative distribution function (cdf) of $Z \sim AbSLG(\alpha,\beta)$ is given by

$$F_Z(z;\alpha,\beta) = \frac{1}{C(\alpha,\beta)(1+e^z)} \begin{pmatrix} e^z(2 + z(\alpha + z^2\beta)(-2 + z\alpha + z^3\beta)) - 2(1+e^z)(\alpha + 3z^2\beta)(-1 + z\alpha + \\ z^3\beta)Log[1+e^z] - 2(1+e^z)((\alpha^2 + 6z(-1+2z\alpha)\beta + 15z^4\beta^2)Li_2[-e^z] + \\ 6\beta((1-4z\alpha - 10z^3\beta)Li_3[-e^z] + (4\alpha + 30z^2\beta)Li_4[-e^z] + \\ 60\beta(-z Li_5[-e^z] + Li_6[-e^z]))) \end{pmatrix} \quad (5)$$

where $Li_n(z) = \sum_{k=1}^{\infty} \frac{z^k}{k^n}$ is the poly-logarithm function (see Prudnikov et al., 1986).

**Proof:** 
$$F_Z(z;\alpha) = \frac{1}{C(\alpha,\beta)} \int_{-\infty}^{z} [(1-\alpha z - \beta z^3)^2 + 1] \frac{e^{-z}}{(1+e^{-z})^2} dz$$

$$= \frac{1}{C(\alpha,\beta)} \left[ \int_{-\infty}^{z} \frac{2e^{-z}}{(1+e^{-z})^2} dz - 2\alpha \int_{-\infty}^{z} \frac{z e^{-z}}{(1+e^{-z})^2} dz + \alpha^2 \int_{-\infty}^{z} \frac{z^2 e^{-z}}{(1+e^{-z})^2} dz - 2\beta \int_{-\infty}^{z} \frac{z^3 e^{-z}}{(1+e^{-z})^2} dz \right.$$

$$\left. + 2\alpha\beta \int_{-\infty}^{z} \frac{z^4 e^{-z}}{(1+e^{-z})^2} dz + \beta^2 \int_{-\infty}^{z} \frac{z^6 e^{-z}}{(1+e^{-z})^2} dz \right]$$

$$= \frac{1}{C(\alpha,\beta)} \left[ 2\left(\frac{e^z}{1+e^z}\right) - 2\alpha\left(\frac{z e^z}{1+e^z} - Log[1+e^z]\right) + \alpha^2 \left\{ z\left(\frac{z e^z}{1+e^z} - 2Log[1+e^z]\right) - 2Li_2[-e^z] \right\} - \right.$$

$$2\beta\left\{ z^2\left(z - \frac{z}{1+e^z} - 3Log[1+e^z]\right) - 6z Li_2[-e^z] + 6 Li_3[-e^z] \right\} +$$

$$2\alpha\beta\left( \frac{z^4 e^z}{1+e^z} - 4z^3 Log[1+e^z] - 12z^2 Li_2[-e^z] + 24z Li_3[-e^z] - 24 Li_4[-e^z] \right) +$$

$$\left. \beta^2 \left\{ z^5 \left(\frac{z e^z}{1+e^z} - 6Log[1+e^z]\right) - 30z^4 Li_2[-e^z] + 120 z(z^2 Li_3[-e^z] - 3z Li_4[-e^z] + 6 Li_5[-e^z]) \right\} - 720 Li_6[-e^z] \right]$$

On simplification we get the result in eqn. (5).

**Corollary 1:** In particular, by taking the limit $\alpha \to \pm\infty$ of $F_Z(z;\alpha,\beta)$ in eqn. (5), we get the cdf of $BLG(2)$ distribution as

$$F_Z(z) = \frac{1}{\pi^2} \left[ 3\left( z\left(\frac{e^z z}{1+e^z} - 2Log[1+e^z]\right) - 2Li_2[-e^z] \right) \right]$$

**Corollary 2:** In particular, by taking the limit $\beta \to \pm\infty$ of $F_Z(z;\alpha,\beta)$ in eqn. (5), we get the cdf of $BLG(6)$ distribution as

$$F_Z(z) = \frac{21}{31(1+e^z)\pi^6} \begin{pmatrix} e^z z^6 - 6z^5 Log[1+e^z](1+e^z) - 30 z^4(1+e^z)Li_2[-e^z] + 120 z^3(1+e^z)Li_3[-e^z] - \\ 360 z^2 Li_4[-e^z](1+e^z) + 720 z Li_5[-e^z](1+e^z) - 720 Li_6[-e^z](1+e^z) \end{pmatrix}$$

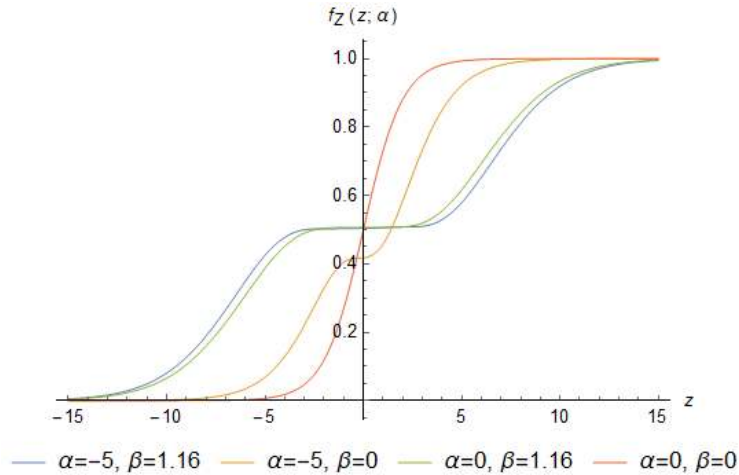

**Figure 2:** The plots of cdf of $AbSLG_2(\alpha,\beta)$ for different values of $\alpha$ and $\beta$

## 3. Moments of AbSLG ($\alpha, \beta$):

**Proposition 3:** The $n^{th}$ moment of $Z \sim AbSLG(\alpha,\beta)$ distribution is given by

$$E(Z^n) = \frac{2}{C(\alpha,\beta)}$$

$$\begin{cases} \left[2\left(1-\frac{1}{2^{n-1}}\right)\Gamma(n+1)\varsigma(n) + \alpha^2\left(1-\frac{1}{2^{n+1}}\right)\Gamma(n+3)\varsigma(n+2) + 2\alpha\beta\left(1-\frac{1}{2^{n+3}}\right)\Gamma(n+5)\varsigma(n+4) + \right. \\ \left. \beta^2\left(1-\frac{1}{2^{n+5}}\right)\Gamma(n+7)\varsigma(n+6)\right]; & \text{if } k \text{ is even} \\ -2\alpha\left(1-\frac{1}{2^n}\right)\Gamma(n+2)\varsigma(n+1) - 2\beta\left(1-\frac{1}{2^{n+2}}\right)\Gamma(n+4)\varsigma(n+3); & \text{if } k \text{ is odd} \end{cases}$$

(6)

*Proof:* Since, the $n^{th}$ moment of the standard logistic distribution (see Balakrishnan, 1992) is given by

$$E(X^n) = \begin{cases} 2\left(1-\frac{1}{2^{n-1}}\right)\Gamma(n+1)\varsigma(n); & \text{when } n \text{ is even} \\ 0; & \text{when } n \text{ is odd} \end{cases} \quad (7)$$

Now, $E(Z^n) = \frac{1}{C(\alpha,\beta)} \int_{-\infty}^{\infty} Z^n \frac{[(1-\alpha z - \beta z^3)^2 + 1]e^{-z}}{(1+e^{-z})^2} dz$

$$= \frac{1}{C(\alpha,\beta)} \left[ \int_{-\infty}^{\infty} \frac{2z^n e^{-z}}{(1+e^{-z})^2} dz - 2\alpha \int_{-\infty}^{\infty} \frac{z^{n+1} e^{-z}}{(1+e^{-z})^2} dz + \alpha^2 \int_{-\infty}^{\infty} \frac{z^{n+2} e^{-z}}{(1+e^{-z})^2} dz - 2\beta \int_{-\infty}^{\infty} \frac{z^{n+3} e^{-z}}{(1+e^{-z})^2} dz \right.$$

$$\left. + 2\alpha\beta \int_{-\infty}^{\infty} \frac{z^{n+4} e^{-z}}{(1+e^{-z})^2} dz + \beta^2 \int_{-\infty}^{\infty} \frac{z^{n+6} e^{-z}}{(1+e^{-z})^2} dz \right]$$

**Case1:** Let $n$ is even, then we get

$$E(Z^n) = \frac{1}{C(\alpha,\beta)} \left[ \int_{-\infty}^{\infty} \frac{2z^n e^{-z}}{(1+e^{-z})^2} dz + \alpha^2 \int_{-\infty}^{\infty} \frac{z^{n+2} e^{-z}}{(1+e^{-z})^2} dz + 2\alpha\beta \int_{-\infty}^{\infty} \frac{z^{n+4} e^{-z}}{(1+e^{-z})^2} dz + \beta^2 \int_{-\infty}^{\infty} \frac{z^{n+6} e^{-z}}{(1+e^{-z})^2} dz \right]$$

$$= \frac{1}{C(\alpha,\beta)} \left[ 4\left(1-\frac{1}{2^{n-1}}\right)\Gamma(n+1)\varsigma(n) + 2\alpha^2\left(1-\frac{1}{2^{n+1}}\right)\Gamma(n+3)\varsigma(n+2) + \right.$$

$$\left. 4\alpha\beta\left(1-\frac{1}{2^{n+3}}\right)\Gamma(n+5)\varsigma(n+4) + 2\beta^2\left(1-\frac{1}{2^{n+5}}\right)\Gamma(n+7)\varsigma(n+6)\right] \quad \text{[Using eqn. (7)]}$$

On simplification we get the result in eqn. (6).

**Case2:** Let $n$ is odd, then we get

$$E(Z^n) = \frac{1}{C(\alpha,\beta)} \left[ -2\alpha \int_{-\infty}^{\infty} \frac{z^{n+1} e^{-z}}{(1+e^{-z})^2} dz - 2\beta \int_{-\infty}^{\infty} \frac{z^{n+3} e^{-z}}{(1+e^{-z})^2} dz \right]$$

$$= \frac{1}{C(\alpha,\beta)} \left[ -4\alpha\left(1-\frac{1}{2^n}\right)\Gamma(n+2)\varsigma(n+1) - 4\beta\left(1-\frac{1}{2^{n+2}}\right)\Gamma(n+5)\varsigma(n+3)\right] \quad \text{[Using eqn. (7)]}$$

On simplification we get the result in eqn. (6).

In particular for $n = 1, 2, 3, 4$ from above, we get

$$E(Z) = -\frac{14(5\pi^2\alpha + 7\pi^4\beta)}{210 + 35\pi^2\alpha^2 + 98\pi^4\alpha\beta + 155\pi^6\beta^2}$$

$$E(Z^2) = \frac{\pi^2(70 + 49\pi^2\alpha^2 + 310\pi^4\alpha\beta + 889\pi^6\beta^2)}{210 + 35\pi^2\alpha^2 + 98\pi^4\alpha\beta + 155\pi^6\beta^2}$$

$$E(Z^3) = -\frac{98\pi^4\alpha + 310\pi^6\beta}{210 + 35\pi^2\alpha^2 + 98\pi^4\alpha\beta + 155\pi^6\beta^2}$$

$$E(Z^4) = \frac{\pi^4(1078 + 1705\pi^2\alpha^2 + 19558\pi^4\alpha\beta + 89425\pi^6\beta^2)}{11(210 + 35\pi^2\alpha^2 + 98\pi^4\alpha\beta + 155\pi^6\beta^2)}$$

$$Var(Z) = \frac{\begin{pmatrix}\pi^2(210 + 35\pi^2\alpha^2 + 98\pi^4\alpha\beta + 155\pi^6\beta^2)(70 + 49\pi^2\alpha^2 + 310\pi^4\alpha\beta + 889\pi^6\beta^2) - \\ 196(5\pi^2\alpha + 7\pi^4\beta)^2\end{pmatrix}}{(210 + 35\pi^2\alpha^2 + 98\pi^4\alpha\beta + 155\pi^6\beta^2)^2}$$

**Corollary 3:** By taking limit $\alpha \to \pm\infty$ in the moments of $AbSLG(\alpha, \beta)$ distribution, we get the moments of $BLG(2)$ distribution as

$$E(Z) \to 0, \quad Var(Z) \to 1.4\pi^2$$

**Corollary 4:** By taking limit $\beta \to \pm\infty$ in the moments of $AbSLG(\alpha, \beta)$ distribution, we get the moments of $BLG(6)$ distribution as

$$E(Z) \to 0, \quad Var(Z) \to 5.7355\pi^2$$

**Proposition 4:** The moment generating function of $Z \sim AbSLG(\alpha, \beta)$ distribution is given by

$$M_Z(t) = \frac{\pi\, Csc[\pi t]}{C(\alpha,\beta)} \begin{pmatrix} -2\alpha + 6\pi^2\beta - t(-2 + (\pi\alpha - \pi^3\beta)^2) + 2\pi((\alpha - \pi^2\beta)(t - \alpha + 3\pi^2\beta)Cot[\pi t] + \\ \pi(-6\beta + t(\alpha - 13\pi^2\beta)(\alpha - 7\pi^2\beta) + 6\pi\beta(t - 4\alpha + 30\pi^2\beta)Cot[\pi t])Csc[\pi t]^2 - \\ 12\pi^3\beta(-2t\alpha + 35\pi^2 t\beta + 30\pi\beta\, Cot[\pi t])Csc[\pi t]^4 + 360\pi^5 t\beta^2 Csc[\pi t]^6 \end{pmatrix}$$

where, $-1 < t < 1$ \hfill (8)

**Proof:** $M_Z(t) = \dfrac{1}{C(\alpha,\beta)} \displaystyle\int_{-\infty}^{\infty} e^{tz}\dfrac{[(1-\alpha z - \beta z^3)^2 + 1]e^{-z}}{(1+e^{-z})^2}dz$

$$= \frac{1}{C(\alpha,\beta)}\left[\int_{-\infty}^{\infty}\frac{2e^{tz-z}}{(1+e^{-z})^2}dz - 2\alpha\int_{-\infty}^{\infty}\frac{ze^{tz-z}}{(1+e^{-z})^2}dz + \alpha^2\int_{-\infty}^{\infty}\frac{z^2 e^{tz-z}}{(1+e^{-z})^2}dz - 2\beta\int_{-\infty}^{\infty}\frac{z^3 e^{tz-z}}{(1+e^{-z})^2}dz \right.$$

$$\left. + 2\alpha\beta\int_{-\infty}^{\infty}\frac{z^4 e^{tz-z}}{(1+e^{-z})^2}dz + \beta^2\int_{-\infty}^{\infty}\frac{z^6 e^{tz-z}}{(1+e^{-z})^2}dz\right]$$

$$= \frac{1}{C(\alpha,\beta)} \left[ 2\{-\pi(-1+t)\,Csc\,(\pi t)\} - 2\alpha\{\pi(-1+\pi(-1+t)\,Cot\,(\pi t))Csc\,(\pi t)\} + \alpha^2 \right.$$

$$\left\{ \frac{1}{2}\pi^2\,Csc\,(\pi t)^3(-\pi(-1+t)(3+Cos\,(2\pi t))+2Sin\,(2\pi t)) \right\} - 2\beta$$

$$\left\{ \frac{1}{4}\pi^3\,Csc\,(\pi t)^4(\pi(-1+t)(23\,Cos\,(\pi t)+Cos\,(3\pi t))-3(5\,Sin\,(\pi t)+Sin\,(3\pi t))) \right\} + 2\alpha\beta$$

$$\left\{ \frac{1}{8}\pi^4\,Csc\,(\pi t)^5(-115\,\pi(-1+t)-\pi(-1+t)(76\,Cos\,(2\pi t)+Cos\,(4\pi t))+88\,Sin\,(2\pi t)+4Sin\,(4\pi t)) \right\} +$$

$$\left. \beta^2 \begin{bmatrix} \frac{1}{32}\pi^6 Csc\,(\pi t)^7(-11774\,\pi(-1+t)-\pi(-1+t)(10543\,Cos\,(2\pi t)+722\,Cos\,(4\pi t)+Cos\,(6\pi t)) \\ +8670\,Sin\,(2\pi t)+1416\,Sin\,(4\pi t)+6Sin\,(6\pi t)) \end{bmatrix} \right]$$

On simplification we get the desired result in eqn. (8).

**Corollary 5:** In particular, by taking the limit $\alpha \to \pm\infty$ of $M_Z(t)$ in eqn. (8), we get the mgf of $BLG(2)$ distribution as

$$M_Z(t) = \frac{3}{2} Csc\,[\pi t]^3 \left( \pi t(3 + Cos\,[2\pi t]) - 2\,Sin\,[2\pi t] \right)$$

**Corollary 6:** In particular, by taking the limit $\beta \to \pm\infty$ of $M_Z(t)$ in eqn. (8), we get the mgf of $BLG(6)$ distribution as

$$M_Z(t) = \frac{21}{992} Csc\,[\pi t]^7 \begin{pmatrix} \pi t(10543\,Cos\,[2\pi t] + 722\,Cos\,[4\pi t] + Cos\,[6\pi t]) - 2(-5887\,\pi t + \\ 4335\,Sin\,[2\pi t] + 708\,Sin\,[4\pi t] + 3\,Sin\,[6\pi t]) \end{pmatrix}$$

**Remark 1:** For putting $t = it$ (where, $i = \sqrt{-1}$) in eqn. (8), we get the characteristic function (cf) of $AbSLG(\alpha,\beta)$ distribution as

$$\phi_Z(t) = \frac{\pi\,Csch\,[\pi t]}{C(\alpha,\beta)} \begin{pmatrix} 2i(\alpha - 3\pi^2\beta) - t(-2 + (\pi\alpha - \pi^3\beta)^2) + 2\pi((-it + \alpha - 3\pi^2\beta)(\alpha - \\ \pi^2\beta)\,Coth\,[\pi t] - \pi(6i\beta + t(\alpha - 13\pi^2\beta)(\alpha - 7\pi^2\beta) - 6\pi\beta(it - 4\alpha + \\ 30\,\pi^2\beta)\,Coth\,[\pi t])\,Csch\,[\pi t]^2 + 12\,\pi^3\beta(2t\alpha - 35\,\pi^2 t\beta + 30\,\pi\beta\,Coth\,[\pi t]) \\ Csch\,[\pi t]^4 - 360\,\pi^5 t\beta^2 Csch\,[\pi t]^6) \end{pmatrix} \quad (9)$$

**Corollary 7:** For taking the limit $\alpha \to +\infty\,or -\infty$ of $\phi_Z(t)$ in eqn. (9), we get the cf of $BLG(2)$ distribution as

$$\phi_Z(t) = -\frac{3}{2} Csch\,[\pi t]^3 \left( \pi t(3 + Cosh\,[2\pi t]) - 2\,Sinh\,[2\pi t] \right)$$

**Corollary 8:** For taking the limit $\beta \to +\infty\,or -\infty$ of $\phi_Z(t)$ in eqn. (9), we get the cf of $BLG(6)$ distribution as

$$\phi_Z(t) = -\frac{21}{992} \left[ Csch\,[\pi t]^7 \begin{pmatrix} \pi t(10543\,Cosh\,[2\pi t] + 722\,Cosh\,[4\pi t] + Cosh\,[6\pi t]) - 2(-5887\,\pi t + \\ 4335\,Sinh\,[2\pi t] + 708\,Sinh\,[4\pi t] + 3\,Sinh\,[6\pi t]) \end{pmatrix} \right]$$

**4. Skewness and Kurtosis of AbSLG $(\alpha, \beta)$:**

**Remark 2:** The expression for skewness and kurtosis of $Z \sim AbSLG(\alpha, \beta)$ distribution are respectively given by

$$\beta_1 = \frac{\left(\begin{array}{l} 5488(5\pi^2\alpha + 7\pi^4\beta)^3 + 2(49\pi^4\alpha + 155\pi^6\beta)(210 + 35\pi^2\alpha^2 + 98\pi^4\alpha\beta + 155\pi^6\beta^2)^2 - \\ 42\pi^2(5\pi^2\alpha + 7\pi^4\beta)(210 + 35\pi^2\alpha^2 + 98\pi^4\alpha\beta + 155\pi^6\beta^2)(70 + 49\pi^2\alpha^2 + \\ 310\pi^4\alpha\beta + 889\pi^6\beta^2) \end{array}\right)^2}{\left(\begin{array}{l} -196(5\pi^2\alpha + 7\pi^4\beta)^2 + \pi^2(210 + 35\pi^2\alpha^2 + 98\pi^4\alpha\beta + 155\pi^6\beta^2)(70 + 49\pi^2\alpha^2 + \\ 310\pi^4\alpha\beta + 889\pi^6\beta^2) \end{array}\right)^3}$$

$$\beta_2 = \frac{\left(\begin{array}{l} -1267728(5\pi^2\alpha + 7\pi^4\beta)^4 - 616(5\pi^2\alpha + 7\pi^4\beta)(98\pi^4\alpha + 310\pi^6\beta)(210 + 35\pi^2\alpha^2 + \\ 98\pi^4\alpha\beta + 155\pi^6\beta^2)^2 + 12936(5\pi^3\alpha + 7\pi^5\beta)^2(210 + 35\pi^2\alpha^2 + 98\pi^4\alpha\beta + 155\pi^6\beta^2) \\ (70 + 49\pi^2\alpha^2 + 310\pi^4\alpha\beta + 889\pi^6\beta^2) + \pi^4(210 + 35\pi^2\alpha^2 + 98\pi^4\alpha\beta + 155\pi^6\beta^2)^3 \\ (1078 + 1705\pi^2\alpha^2 + 19558\pi^4\alpha\beta + 89425\pi^6\beta^2) \end{array}\right)}{11\left(\begin{array}{l} -196(5\pi^2\alpha + 7\pi^4\beta)^2 + \pi^2(210 + 35\pi^2\alpha^2 + 98\pi^4\alpha\beta + 155\pi^6\beta^2)(70 + 49\pi^2\alpha^2 + \\ 310\pi^4\alpha\beta + 889\pi^6\beta^2) \end{array}\right)^2}$$

**Corollary 9:** For taking limit $\alpha \to \pm\infty$ in the results of $AbSLG(\alpha, \beta)$ distribution, we get the skewness and kurtosis of $BLG(2)$ distribution as

$$\beta_1 \to 0, \quad \beta_2 \to 2.2595$$

**Corollary 10:** For taking limit $\beta \to \pm\infty$ in the results of $AbSLG(\alpha, \beta)$ distribution, we get the skewness and kurtosis of $BLG(6)$ distribution as

$$\beta_1 \to 0, \quad \beta_2 \to 1.5944$$

**Remark 3:** The density function (4), of model $AbSLG(\alpha, \beta)$ can be represented as sum of two functions, as shown below

$$f_Z(z; \alpha, \beta) = \frac{(2 + \alpha^2 z^2 + 2z^4\alpha\beta + z^6\beta^2)}{C(\alpha, \beta)} \frac{e^{-z}}{(1 + e^{-z})^2} - \frac{(2\alpha z + 2z^3\beta)}{C(\alpha, \beta)} \frac{e^{-z}}{(1 + e^{-z})^2} \quad (10)$$

In equation (10) the 1$^{st}$ part is symmetric and the second part is asymmetric one and the symmetric part, which is defined below, is symbolically denoted by $SCAbSLG(\alpha, \beta)$.

**Definition 2:** A random variable $Z \sim SCAbSLG(\alpha, \beta)$ distribution, if its density function given by

$$f_1(z; \alpha, \beta) = \frac{(2 + \alpha^2 z^2 + 2z^4\alpha\beta + z^6\beta^2)}{C(\alpha, \beta)} \frac{e^{-z}}{(1 + e^{-z})^2} \quad (11)$$

where, $\alpha, \beta \in R$, then $Z$ is a symmetric-component random variable of the model $AbSLG(\alpha, \beta)$. We denote it as $Z \sim SCAbSLG(\alpha, \beta)$

**Some properties of SCAbSLG $(\alpha, \beta)$:**

(i) If $\alpha = \beta = 0$, then we get the standard logistic distribution.

(ii) If $\alpha, \beta \to \pm\infty$, then $Z \xrightarrow{d} BLG$

**Remark 4:** The density function (11) is a mixture between a logistic density and a bimodal-logistic density, as given below

$$f_1(z;\alpha,\beta) = \frac{2}{C(\alpha,\beta)}\frac{e^{-z}}{(1+e^{-z})^2} + \frac{(\alpha^2 z^2 + 2z^4 \alpha\beta + z^6 \beta^2)}{C(\alpha,\beta)}\frac{e^{-z}}{(1+e^{-z})^2}$$

**Proposition 5:** The cdf of $Z \sim AbSLG(\alpha,\beta)$ distribution is given by

$$F_1(z;\alpha,\beta) = \frac{1}{C(\alpha,\beta)(1+e^z)}\begin{pmatrix} e^z(2+z^2(\alpha+z^2\beta)^2) - 2(1+e^z)z(\alpha+z^2\beta)(\alpha+3z^2\beta)Log[1+e^z] - 2(1+e^z)((\alpha^2 + \\ 12z^2\alpha\beta + 15z^4\beta^2)Li_2[-e^z] + 12\beta(-(2z\alpha+5z^3\beta)Li_3[-e^z] + (2\alpha+15z^2\beta)Li_4[-e^z] + \\ 30\beta(-zLi_5[-e^z]+Li_6[-e^z]))) \end{pmatrix}$$

(12)

*Proof:* $F_Z(z;\alpha) = \frac{1}{C(\alpha,\beta)}\int_{-\infty}^{z}(2+\alpha^2 z^2 + 2z^4\alpha\beta + z^6\beta^2)\frac{e^{-z}}{(1+e^{-z})^2}dz$

$= \frac{1}{C(\alpha,\beta)}\left[\int_{-\infty}^{z}\frac{2e^{-z}}{(1+e^{-z})^2}dz + \alpha^2\int_{-\infty}^{z}\frac{z^2 e^{-z}}{(1+e^{-z})^2}dz + 2\alpha\beta\int_{-\infty}^{z}\frac{z^4 e^{-z}}{(1+e^{-z})^2}dz + \beta^2\int_{-\infty}^{z}\frac{z^6 e^{-z}}{(1+e^{-z})^2}dz\right]$

$= \frac{1}{C(\alpha,\beta)}\left[2\left(\frac{e^z}{1+e^z}\right) + \alpha^2\left\{z\left(\frac{ze^z}{1+e^z} - 2Log[1+e^z]\right) - 2Li_2[-e^z]\right\} + \right.$

$2\alpha\beta\left(\frac{z^4 e^z}{1+e^z} - 4z^3 Log[1+e^z] - 12z^2 Li_2[-e^z] + 24z Li_3[-e^z] - 24 Li_4[-e^z]\right) +$

$\left.\beta^2\left\{z^5\left(\frac{ze^z}{1+e^z} - 6Log[1+e^z]\right) - 30z^4 Li_2[-e^z] + 120z(z^2 Li_3[-e^z] - 3z Li_4[-e^z] + 6 Li_5[-e^z]) - 720 Li_6[-e^z]\right\}\right]$

On simplification we get the result in eqn. (12).

**Proposition 6:** The mgf of $Z \sim AbSLG(\alpha,\beta)$ distribution is given by

$$M_Z(t) = \frac{\pi\, Csc[\pi t]}{C(\alpha,\beta)}\begin{pmatrix} -t(-2+(\pi\alpha-\pi^3\beta)^2) + 2\pi(-(\alpha-3\pi^2\beta)(\alpha-\pi^2\beta)Cot[\pi t] + \pi(t(\alpha-13\pi^2\beta)(\alpha - \\ 7\pi^2\beta) + 12\pi\beta(-2\alpha+15\pi^2\beta)Cot[\pi t])Csc[\pi t]^2 - 12\pi^3\beta(-2t\alpha+35\pi^2 t\beta + \\ 30\pi\beta\, Cot[\pi t])Csc[\pi t]^4 + 360\,\pi^5 t\beta^2 Csc[\pi t]^6) \end{pmatrix}$$

(13)

*Proof:* $M_Z(t) = \frac{1}{C(\alpha,\beta)}\int_{-\infty}^{\infty}e^{tz}\frac{(2+\alpha^2 z^2 + 2z^4\alpha\beta + z^6\beta^2)e^{-z}}{(1+e^{-z})^2}dz$

$$= \frac{1}{C(\alpha,\beta)} \left[ \int_{-\infty}^{\infty} \frac{2e^{tz-z}}{(1+e^{-z})^2} dz + \alpha^2 \int_{-\infty}^{\infty} \frac{z^2 e^{tz-z}}{(1+e^{-z})^2} dz + 2\alpha\beta \int_{-\infty}^{\infty} \frac{z^4 e^{tz-z}}{(1+e^{-z})^2} dz + \beta^2 \int_{-\infty}^{\infty} \frac{z^6 e^{tz-z}}{(1+e^{-z})^2} dz \right]$$

$$= \frac{1}{C(\alpha,\beta)} \left[ 2\{-\pi(-1+t)\,Csc\,(\pi t)\} + \alpha^2 \left\{ \frac{1}{2}\pi^2\,Csc\,(\pi t)^3(-\pi(-1+t)(3+Cos\,(2\pi t))+2\,Sin\,(2\pi t))\right\} + \right.$$

$$2\alpha\beta \left\{ \frac{1}{8}\pi^4\,Csc\,(\pi t)^5(-115\pi(-1+t)-\pi(-1+t)(76\,Cos\,(2\pi t)+Cos\,(4\pi t))+88\,Sin\,(2\pi t)+4\,Sin\,(4\pi t))\right\} +$$

$$\left. \beta^2 \left\{ \frac{1}{32}\pi^6 Csc\,(\pi t)^7 \begin{pmatrix} (-11774\pi(-1+t)-\pi(-1+t)(10543\,Cos\,(2\pi t)+722\,Cos\,(4\pi t)+Cos\,(6\pi t)) \\ +8670\,Sin\,(2\pi t)+1416\,Sin\,(4\pi t)+6\,Sin\,(6\pi t)) \end{pmatrix} \right\} \right]$$

On simplification we get the desired result in eqn. (13).

**Algorithm:**

The acceptance-rejection algorithm can be used to simulate random numbers from $AbSLG_2(\alpha,\beta)$ distribution.

Let $g(x;\alpha,\beta)$, the density function of $X \sim AbSLG_2(\alpha,\beta)$ and $g_1(x;\alpha,\beta)$, the density function of $Y \sim SCAbSLG_2(\alpha,\beta)$, with

$$S = \underset{x}{Sup}\,\frac{g(x;\alpha,\beta)}{g_1(x;\alpha,\beta)} = \frac{2+\sqrt{2}}{2}$$

Then a random variable from $AbSLG_2(\alpha,\beta)$ distribution can be generated in the following steps:

(i) Generate $Y \sim SCAbSLG_2(\alpha,\beta)$ and $U \sim U(0,1)$, and they are independent.

(ii) If $U < \frac{1}{S}\frac{g(x;\alpha,\beta)}{g_1(x;\alpha,\beta)}$, and set $X=Y$; otherwise go back to step(i) and continue the process.

### 5. Estimation and Applications of AbSLG $(\alpha,\beta)$:

We present the problem of parameter estimation of a location and scale extension of $AbSLG(\alpha,\beta)$ distribution. If $Z \sim AbSLG(\alpha,\beta)$ then $Y = \mu + \sigma Z$ is said to be the location ($\mu$) and scale ($\sigma$) extension of $Z$ and has the pdf is given by

The location ($\mu$) and scale ($\sigma$) extension of $Y = \mu + \sigma Z$ for $Z \sim AbSLG(\alpha,\beta)$, $\mu \in R$ and $\sigma > 0$ has the following pdf given by

$$f_Y(y;\mu,\sigma,\alpha,\beta) = \frac{\left(\left[1-\alpha\left(\frac{y-\mu}{\sigma}\right)-\beta\left(\frac{y-\mu}{\sigma}\right)^3\right]^2+1\right)}{C(\alpha,\beta)} \frac{e^{-\left(\frac{y-\mu}{\sigma}\right)}}{\sigma\left[1+e^{-\left(\frac{y-\mu}{\sigma}\right)}\right]^2}; y,\mu,\alpha,\beta \in R; \sigma > 0 \quad (14)$$

We denote it by, $Y \sim AbSLG(\mu, \sigma, \alpha, \beta)$. Then, the log-likelihood function is defined by

$$l(\theta) = \sum_{i=1}^{n} \log\left[\left\{1 - \alpha\left(\frac{y_i - \mu}{\sigma}\right) - \beta\left(\frac{y_i - \mu}{\sigma}\right)^3\right\}^2 + 1\right] - n\log(210 + 35\pi^2\alpha^2 + 98\pi^4\alpha\beta + 155\pi^6\beta^2) + n\log[105]$$

$$- n\log(\sigma) - \sum_{i=1}^{n} \frac{(y_i - \mu)}{\sigma} - 2\sum_{i=1}^{n} \log\left[1 + Exp\left[\frac{-(y_i - \mu)}{\sigma}\right]\right] \quad (15)$$

For analysis purposes, we consider three datasets. The first one is obtained from the website http://users.stat.umn.edu/sandy/courses/8061/datasets/lakes.lsp which is related to N latitude degrees in 69 samples from world lakes, which appear in Column 5 of the Diversity data set. The second one is acquired from the website http://www.globalfindata.com which is the exchange rate data of the United Kingdom Pound to the United States Dollar from 1800 to 2003. The third one consists of the velocities of 82 distant galaxies, diverging from our own galaxy in the website http://www.stats.bris.ac.uk/~peter/mixdata.

For the above datasets we will fit the following models: the proposed model $AbSLG(\alpha, \beta, \mu, \sigma)$, the normal model $N(\mu, \sigma^2)$, the logistic model $LG(\mu, \beta)$, the Laplace model $La(\mu, \beta)$, the skew-normal model $SN(\lambda, \mu, \sigma)$ of Azzalini (1985), the skew-logistic model $SLG(\lambda, \mu, \beta)$ of Wahed and Ali (2001), the skew-Laplace distribution $SLa(\lambda, \mu, \beta)$ of Nekoukhou and Alamatsaz (2012), the alpha-skew-normal model $ASN(\alpha, \mu, \sigma)$ of Elal-Olivero (2010), the alpha-skew-Laplace model $ASLa(\alpha, \mu, \beta)$ of Harandi and Alamatsaz (2013), the alpha-skew-logistic model $ASLG(\alpha, \mu, \beta)$ of Hazarika and Chakraborty (2014), and beta skew logistic model $bSLG(\beta, \mu, \sigma)$.

Employing the maximum likelihood procedure we can estimate the parameters of the proposed $AbSLG(\alpha, \beta, \mu, \sigma)$ model. AIC and BIC are used for model comparison.

**Table 1:** MLE's, log-likelihood, AIC and BIC for N latitude degrees in 69 samples from world lakes.

| Distributions | $\mu$ | $\sigma$ | $\lambda$ | $\alpha$ | $\beta$ | $\log L$ | AIC | BIC |
|---|---|---|---|---|---|---|---|---|
| $N(\mu, \sigma^2)$ | 45.165 | 9.549 | -- | -- | -- | -253.599 | 511.198 | 515.666 |
| $LG(\mu, \beta)$ | 43.639 | -- | -- | -- | 4.493 | -246.645 | 497.29 | 501.758 |
| $SN(\lambda, \mu, \sigma)$ | 35.344 | 13.698 | 3.687 | -- | -- | -243.036 | 492.072 | 498.774 |
| $SLG(\lambda, \mu, \beta)$ | 36.787 | -- | 2.828 | -- | 6.417 | -239.053 | 490.808 | 490.808 |
| $La(\mu, \beta)$ | 43.00 | -- | -- | -- | 5.895 | -239.248 | 482.496 | 486.964 |
| $ASLG(\alpha, \mu, \beta)$ | 49.087 | -- | -- | 0.861 | 3.449 | -237.351 | 480.702 | 487.404 |
| $SLa(\lambda, \mu, \beta)$ | 42.30 | -- | 0.255 | -- | 5.943 | -236.900 | 479.799 | 486.501 |
| $ASLa(\alpha, \mu, \beta)$ | 42.30 | -- | -- | -0.220 | 5.439 | -236.079 | 478.159 | 484.861 |
| $ASN(\alpha, \mu, \sigma)$ | 52.147 | 7.714 | -- | 2.042 | -- | -235.370 | 476.739 | 483.441 |
| $bSLG(\beta, \mu, \sigma)$ | 42.285 | 3.172 | -- | -- | -0.012 | -235.199 | 476.398 | 483.100 |
| $AbSLG(\alpha, \beta, \mu, \sigma)$ | 45.487 | 2.755 | -- | 0.737 | -0.029 | **-230.710** | **469.416** | **478.352** |

**Table 2:** MLE's, log-likelihood, AIC and BIC for the exchange rate data of the United Kingdom Pound to the United States Dollar from 1800 to 2003.

| Distributions | $\mu$ | $\sigma$ | $\lambda$ | $\alpha$ | $\beta$ | $\log L$ | AIC | BIC |
|---|---|---|---|---|---|---|---|---|
| $N(\mu,\sigma^2)$ | 4.117 | 1.381 | -- | -- | -- | -355.265 | 714.529 | 721.165 |
| $SN(\lambda,\mu,\sigma)$ | 3.589 | 1.478 | 0.501 | -- | -- | -355.217 | 716.434 | 726.388 |
| $LG(\mu,\beta)$ | 4.251 | -- | -- | -- | 0.753 | -351.192 | 706.385 | 713.021 |
| $SLG(\lambda,\mu,\beta)$ | 5.36 | -- | -2.371 | -- | 1.018 | -341.391 | 688.782 | 698.736 |
| $La(\mu,\beta)$ | 4.754 | -- | -- | -- | 0.971 | -339.315 | 682.63 | 689.265 |
| $ASN(\alpha,\mu,\sigma)$ | 3.656 | 0.883 | -- | -3.504 | -- | -317.946 | 641.892 | 651.847 |
| $SLa(\lambda,\mu,\beta)$ | 4.855 | -- | 1.506 | -- | 1 | -311.318 | 628.636 | 638.59 |
| $ASLG(\alpha,\mu,\beta)$ | 3.764 | -- | -- | -2.025 | 0.403 | -301.963 | 609.927 | 619.881 |
| $bSLG(\beta,\mu,\sigma)$ | 4.723 | 0.32 | -- | -- | 0.024 | -301.683 | 609.366 | 619.32 |
| $ASLa(\alpha,\mu,\beta)$ | 4.861 | -- | -- | 0.539 | 0.677 | -301.443 | 608.885 | 618.84 |
| $AbSLG(\alpha,\beta,\mu,\sigma)$ | 4.447 | 0.248 | -- | -1.465 | 0.069 | -264.235 | **536.47** | **549.743** |

**Table 3:** MLE's, log-likelihood, AIC and BIC for the velocities of 82 distant galaxies, diverging from our own galaxy.

| Distributions | $\mu$ | $\sigma$ | $\lambda$ | $\alpha$ | $\beta$ | $\log L$ | AIC | BIC |
|---|---|---|---|---|---|---|---|---|
| $N(\mu,\sigma^2)$ | 20.832 | 4.54 | -- | -- | -- | -240.417 | 484.833 | 489.646 |
| $SN(\lambda,\mu,\sigma)$ | 24.61 | 5.907 | -1.395 | -- | -- | -239.21 | 484.42 | 491.64 |
| $SLG(\lambda,\mu,\beta)$ | 21.532 | -- | -0.154 | -- | 2.219 | -233.314 | 472.628 | 479.849 |
| $LG(\mu,\beta)$ | 21.075 | -- | -- | -- | 2.204 | -233.649 | 471.299 | 476.113 |
| $ASN(\alpha,\mu,\sigma)$ | 17.417 | 3.869 | -- | -1.656 | -- | -230.088 | 466.175 | 473.395 |
| $SLa(\lambda,\mu,\beta)$ | 20.846 | -- | 1.002 | -- | 2.997 | -228.829 | 463.658 | 470.878 |
| $La(\mu,\beta)$ | 20.838 | -- | -- | -- | 2.997 | -228.83 | 461.66 | 466.474 |
| $ASLG(\alpha,\mu,\beta)$ | 18.482 | -- | -- | -0.833 | 1.646 | -224.877 | 455.754 | 462.974 |
| $ASLa(\alpha,\mu,\beta)$ | 19.473 | -- | -- | -0.842 | 1.805 | -220.789 | 447.578 | 454.798 |
| $bSLG(\beta,\mu,\sigma)$ | 20.779 | 1.279 | -- | -- | -0.017 | -220.747 | 447.494 | 454.714 |
| $AbSLG(\alpha,\beta,\mu,\sigma)$ | 20.317 | 1.432 | -- | -0.454 | 0.022 | -218.802 | **445.604** | **455.231** |

From the Tables 1 - 3, it is seen that the proposed alpha-beta-skew-logistic $AbSLG(\alpha,\beta,\mu,\sigma)$ distribution provides best fit to the data set under consideration in terms of all the criteria, namely the log-likelihood, the AIC as well as the BIC. From the figures 3 – 5 below, which are the plots of observed (in histogram) and expected densities (lines), also confirms our finding.

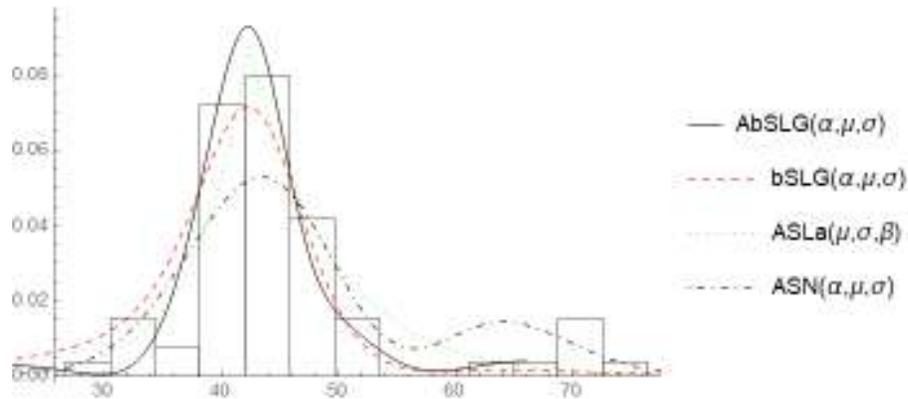

**Figure 3:** Plots of observed and expected densities of some distribution for N latitude degrees in 69 samples from world lakes.

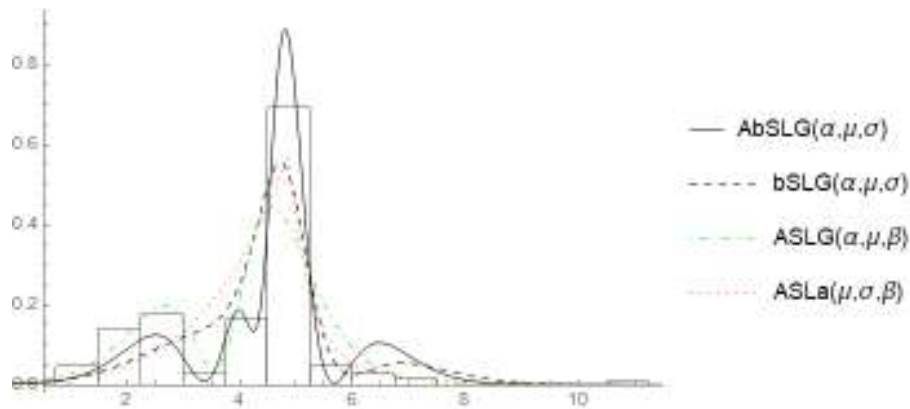

**Figure 4:** Plots of observed and expected densities of some distributions for the exchange rate data of the United Kingdom Pound to the United States Dollar from 1800 to 2003.

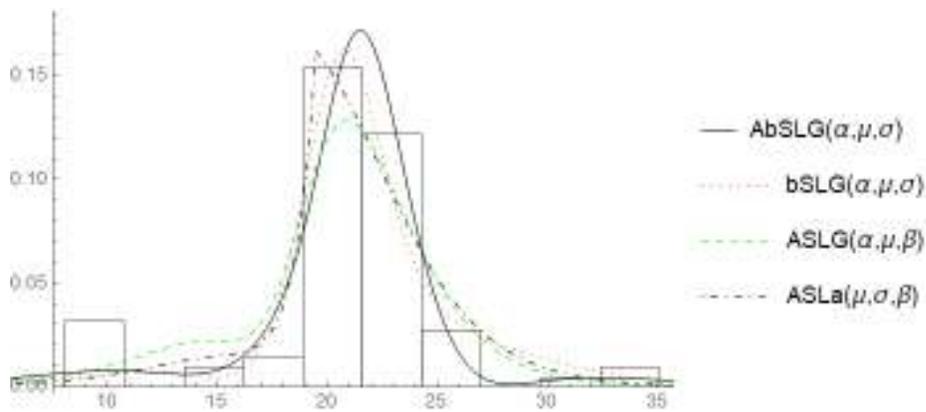

**Figure 5:** Plots of observed and expected densities of some distributions for the velocities of 82 distant galaxies, diverging from our own galaxy.

**Likelihood Ratio Test:**

Further, the likelihood ratio (LR) test is carried to discriminate between the nested models: $LG(\mu,\beta)$ and $AbSLG(\alpha,\beta,\mu,\sigma)$, under the Null hypothesis, $H_0: \alpha = \beta = 0$, that is the sample is drawn from $LG(\mu,\beta)$; against the alternative $H_1: \alpha = \beta \neq 0$, that is the sample is drawn from $AbSLG(\alpha,\beta,\mu,\sigma)$.

**Table 3:** The values of LR test statistic.

|  | Dataset 1 | Dataset 2 | Dataset 3 |
|---|---|---|---|
| LR test value | 31.7894 | 173.914 | 29.694 |

Since, in all the four datasets, the values of LR test statistic (i.e., 31.7894, 173.914, and 29.694 respectively) exceed the 99% critical value (i.e., 6.635). Thus, it favours the alternative hypothesis that the sampled data comes from $AbSLG(\alpha,\beta,\mu,\sigma)$, not from $LG(\mu,\beta)$.

## 6. Conclusion

In this study a new skew distribution is proposed and is known as alpha beta skew logistic distribution which includes unimodal as well as bimodal behavior and some of its basic properties are investigated. From the computation, it is examined that the proposed $BAbSLG_2(\alpha,\mu,\beta)$ distribution provides best fitting to the data set under consideration in terms of all the criteria, namely log-likelihood, AIC and BIC. The plots of observed and expected densities presented above, also confirms our findings.